# Note on the *ABC* Conjecture
N. A. Carella, February 2006

*Abstract:* This note imparts heuristic arguments and theoretical evidences that contradict the *abc* conjecture over the rational integers. In addition, the rudimentary details for transforming this problem into the domain of equidistribution theory are provided.

**Introduction**
The *abc* conjecture proclaims a link between the additive and the multiplicative structures of the integer solutions of the equation Diophantine $x + y = z$. The precise statement is as follows.

***Conjecture* 1.** (Masser and Oesterle, 1985) Let $a, b, c \in \mathbb{N}$ be integers. Given a real number $\varepsilon > 0$, there exists a constant $c_\varepsilon$ such that if $a + b = c$, $\gcd(a, b, c) = 1$, then

$$\max\{|a|,|b|,|c|\} \leq c_\varepsilon \prod_{p|abc} p^{1+\varepsilon}. \qquad (1)$$

Weaker versions, and generalizations over numbers fields, and functions fields, et cetera and related ideas are considered in [BA04], [GE02], [BB94], [VJ00], [LG90] and the references therein.

Let $n \in \mathbb{N}$ be an integer. The radical and the index are defined by the expressions $rad(n) = \prod_{p|n} p$, and $\gamma(a,b) = \dfrac{\max\{\log|a|,|b|,|c|\}}{\log rad(abc)}$ respectively. In term of limit inequality (1) states that the conjectured maximal limit of the set of all indices is given by

$$\limsup_{\gcd(a,b)=1} \gamma(a,b) = 1. \qquad (2)$$

The limit aspects of the abc conjecture and associated technique are explored in [BR00], [FA98], [BF95], [GN99]. The later proves that every point in the interval [1/3, 36/37] is a limit point of the index. Now since the index $\gamma(a, b) \in [1/3, A] \subset [1/3, \infty]$ for some $A \geq 1$, and it has an average value of 1, see [DK05], it should have a distribution on the half line [1/3, ∞] and somewhat centered on 1, and not strictly confined to the interval [1/3, 1].

The heuristic/theoretical evidences herein do not contradict a weaker version such as $\limsup_{\gcd(a,b)=1} \gamma(a,b) \leq c_0$, where $c_0 > 1$ is some constant. This is not a severe impediment, since even a weak abc result is a powerful tool in number theory. Extensive analysis and applications are given in [BI06].

**Heuristic/Theoretical Evidences**

The construction of triples $a$, $b$, $c$ by continued fraction approximations of the algebraic integers $a^{1/n}$, $a$, $n \geq 2$, [BB94], [NJ03], and integer relations $Ax_1 + Bx_2 + Cx_3$, with $A$, $B$, $C$ sufficiently smooth [DR04] provide quite general methods of generating triples of large indices. In this section a different technique of generating triples $a$, $b$, $c$ of radical indices $\gamma(a, b) > 1$ will be considered. This is motivated by certain entries of the form $N = q + 1 - t$ in the table [NJ06] of large indices $\gamma(a, b) > 1$, and the fact that there are groups of points of algebraic curves of highly smooth orders $N$.

For example, the record abc triple $(2, 109 \cdot 3^{10}, 23^5)$ of index $\gamma(a, b) = 1.62991$ was originally obtained from the continued fraction approximation 23/9 of $109^{1/5} = [2; 1, 1, 4, 77733, …]$. The triple can also be viewed as parameters of an algebraic curve of genus $g = 1$. Specifically, $23^5 + 1 - 3 = 3^{10} \cdot 109$ arises as the order $N = \#E(\mathbf{F}_q)$ of the group $E(\mathbf{F}_q)$ of some elliptic curve $E : y^2 = x^3 + ax + b$ over $\mathbf{F}_q$, $q = 23$ or $23^5$.

Perhaps, this is not surprising since a lot of the abc related ideas spring from the theory of elliptic curves. For example, Szapiro discriminant conjecture, and the size of the Shafarevich-Tate groups. These are

$$\Delta \ll N^{6+\varepsilon}, \text{ and } \#\text{Ш} \ll N^{1/2+\varepsilon} \tag{3}$$

respectively, see [GS95], [WR98].

Let $\bigcup_q I_q$ be the union of the Hasse intervals $I_q = [q + 1 - 2g\sqrt{q}, q + 1 + 2g\sqrt{q}]$ of algebraic curves over $\mathbf{F}_q$ of genus $g \geq 1$, and let $d_n = \max \{ p_{k+1} - p_k : k \leq n \}$ be the prime gap between the $n$th primes $p_n$ and $p_{n+1}$. It is obvious that the equality $\bigcup_q I_q = \mathbb{N}$ occurs if the prime gap between primes $p_n$ and $p_{n+1}$ is bounded by $4g\sqrt{p_n}$. The current best estimates of the prime gap are $2\sqrt{p_n}\log(p_n)$ modulo the Riemann hypothesis, and $d_n \leq p_n^{.525}$ unconditionally, [BA01]. Accordingly, the union of the Hasse intervals of algebraic curves of genus $d_n \leq p_n^{.525}$ is an unconditional cover the positive integers line, e.g., $\bigcup_q I_q = \mathbb{N}$, see [BS04] for application and discussion.

The smoothness and primality of integers in short intervals $[X, X + Y]$, and Hasse intervals of small genus $g = 1, 2$ has been extensively investigated in the last two decades, see [HN91], [HN99], [SC03], [LN87], LN02]. These works show that there are sufficiently many smooth integers in these short intervals.

The application undertaken here calls for vastly weaker condition, it assumes Hasse intervals of algebraic curves up to exponentially large genera as required. Thus the intervals $I_q = [\, q+1-2g\sqrt{q},\, q+1+2g\sqrt{q}\,]$ are almost certain to support algebraic curves of highly smooth orders (and very few prime factors) for infinitely many prime powers.

Let **F** be a finite field of $q = p^n$ elements, and consider an algebraic curve $C : y^2 = f(x)$ of genus $g > 0$. Further, let $N = \#C(\mathbf{F}) = q + 1 - t$ be the cardinality of its set of **F**-rational points, and put

$$a = N, \quad b = t - 1, \quad c = a + b. \tag{4}$$

***Heuristic Claim 2.*** Let $\varepsilon < 1/12$. Then the abc inequality (1) fails for infinitely many rational integer triples $a, b, c$ with $\gcd(a, b, c) = 1$.

***Sketch of the argument.*** Let $g \le q^\delta$, $\delta < 1/12$, and suppose that $N = \#C(\mathbf{F}_q) = q + 1 - t$. By Theorem 7 the interval $[\, q+1-2g\sqrt{q},\, q+1+2g\sqrt{q}\,]$ contains significantly more than $q^{1/2}$ smooth integers, and by Theorem 10, there are algebraic curves of genus $g = 2$ and order $N$ in this internal for almost every $N$. Ergo it is reasonable to expect smooth orders $N$ such that $\mathrm{rad}(N) \le q^{1/3}$. Under these conditions, it follows that

$$q \le c_\varepsilon \prod_{\text{prime } r \mid abc} r^{1+\varepsilon} \le c_\varepsilon \bigl(\mathrm{rad}(N(t-1)q)\bigr)^{1+\varepsilon} \le c_\varepsilon \bigl(p^{1+(\delta+5/6)n}\bigr)^{1+\varepsilon} \le c_\varepsilon \bigl(p^{1+11n/12}\bigr)^{1+\varepsilon}. \tag{5}$$

Taking logarithm and simplifying yield

$$\frac{12}{11+12/n} \le 1 + \varepsilon + \frac{12 \log c_\varepsilon}{(11+12/n) n \log p}. \tag{6}$$

As $n, p \to \infty$ the claim follows. ∎

There is a well known algorithm for constructing algebraic curves of specified orders $N$ such that $4q = x^2 + dy^2$, and $N = (x-1)^2 + dy^2$, some $d > 0$. This algorithm is very successful in constructing elliptic curves of prime and square-free orders, (the actual construction of the curves is not required). This same algorithm should be successful in finding algebraic curves of $k$-powerful or nearly $k$-powerful orders. The interval $[1, X]$ contains approximately $cX^{1/k} + O(X^{1/(k+1)})$ $k$-powerful integers, [IV88], but the current estimates for the number of $k$-powerful numbers in short intervals $[\, X - cX^{1/2+\varepsilon},\, X + cX^{1/2+\varepsilon}\,]$ are very weak or nonexistent, see the Appendix. However, it should be an increasing function of $X$. Moreover, large indices should be possible.

**Equidistribution Approach**
An investigation of the smoothness and statistical property of the integers $N = q + 1 - a_n$ as a function of $q$ and the Fourier coefficients $a_n$ of the automorphic $L$-function

$$L(s) = \sum_{n=1}^{\infty} \frac{a_n}{n^s} = \prod_{\gcd(p,N_0)>1}(1-a_p p^{-s})^{-1} \prod_{\gcd(p,N_0)=1}(1-a_p p^{-s} + p^{1-2s})^{-1}, \qquad (7)$$

where $N_0$ is the conductor of the elliptic curve could yield surprising results in connection with the abc conjecture. A related but simpler approach will be explored next.

To obtain the elementary details of transforming the abc conjecture into an equidistribution problem, [IK04, Chapter 21], let $k = 2, 3$ or $4$, and let $\| x \|$ be the nearest integer to $x > 0$. It is a simple observation that every integer $n \in \mathbb{N}$ has a $k$-representation

$$n = \| n^{1/k} \|^k + b_n, \ 0 \leq | b_n | < 2kn^{(k-1)/k}. \qquad (8)$$

Several well known families of equations can be viewed as special cases of this construction. Some of these are
(i) The Ramanujan-Nagell equation $x^2 + d = p_1^{e_1} p_2^{e_{21}} \cdots p_v^{e_v}$, this have been solved for various parameters, see the literature.
(ii) The square-cubic difference equation $x^3 = y^2 + d$, the relevance here is the Hall's conjecture, which states that $|x^3 - y^2| > cx^{1/2}$. But this has been contradicted in [DV82] who showed that $|x^3 - y^2| < cx^{1/2}$ infinitely often. Moreover it is not implausible to have infinitely many very small or very smooth $d$.
(iii) The Gaussian integers $n = a^2 + b^2$. It has an increasing number of representations as a sum of two squares depending on the prime factors in $n$, and can be recursively generated, i. e., $5 = 2^2 + 1^2$, $5^2 = (2 \cdot 1 \pm 1 \cdot 1)^2 + (2 \cdot 1 \mp 1 \cdot 1)^2$, $5^3 = (2 \cdot 3 \pm 1 \cdot 4)^2 + (4 \cdot 3 \mp 2 \cdot 1)^2$, ... .

The literature on these families of equations is vast, but there is nothing from the equidistribution point of view.

For $k = 2$, the associated Dirichlet series

$$L(s) = \sum_{n=1}^{\infty} \frac{b_n}{n^s} = 1 + \frac{1}{2^s} - \frac{1}{3^s} + \frac{1}{5^s} + \frac{2}{6^s} - \frac{2}{7^s} - \frac{1}{8^s} + \frac{1}{10^s} + \frac{2}{11^s} + \frac{3}{12^s} \pm \frac{4}{13^s} - \frac{3}{14^s} + \cdots \qquad (9)$$

is absolutely convergent in the complex half plane $\text{Re}(s) > 3/2$. This series probably does not have an Euler product

$$L(s) = \sum_{n=1}^{\infty} \frac{b_n}{n^s} \stackrel{?}{=} \prod_{\text{prime } p}(1-\psi(p)p^{-s})^{-1} \qquad (10)$$

(for some arithmetic function $\psi$) since the Hecke relations

(i) $b_p b_q$     $\gcd(p, q) = 1$,     (ii) $b_{p^n} = b_p b_{p^{n-1}} - b_{p^{n-2}}$,

do not hold, see [BK00] for related work. However, it still has very nice properties.

**Theorem 3.** For all large $x > 0$, and $n \leq x$, the integers $b_n$ are uniformly distributed with respect to some measure in the interval $[-x^{1/2}, x^{1/2}]$.

**Theorem 4.** Let $\delta > 0$ be a sufficiently small number, and let $q = \| q^{1/2} \|^2 + b_q$ be the representations of the sequence of prime powers $q = p^{n+3}$, $n \geq 0$, $p \geq 2$. Suppose that $0 < |b_q| \leq 2q^{1/2-\delta}$ holds for infinitely many $q$. Then the abc conjecture fails.

Proof: The appropriate radical is $\mathrm{rad}(q \cdot \| q^{1/2} \| \cdot b_q) \leq p^{1 + (n+3)/2 + (1/2 - \delta)(n+3)} = p^{1 + (1 - \delta)(n+3)}$. Hence

$$\lim_{n,p \to \infty} \frac{\log q}{\log \mathrm{rad}(q \cdot \| q^{1/2} \| \cdot b_q)} \geq \lim_{n,p \to \infty} \frac{(n+3)\log p}{((1-\delta)(n+3)+1)\log(p)} = \frac{1}{1-\delta} > 1. \quad (11)$$

This proves the contention. ∎

Does the uniform distribution of the integers $b_n$ with $n \leq x$ in the interval $[-x^{1/2}, x^{1/2}]$ imply the uniform distribution of the integers $b_q$ with $q = p^{n+3} \leq x$ in the interval $[-x^{1/2}, x^{1/2}]$. A proof that the sequence $\{ b_q \}$ is uniform distributed immediately disproves the abc conjecture. But the full strength of uniform distribution of the values $b_q$ with respect to some measure is not required. The above hypothesis uses a weaker condition on the distribution of the values $b_q$. For example, the existence of the set

$$\#\{ p^{n+3} \leq x : 0 < b_{p,n} \leq p^{(n+3)(1/2-\delta)} \} \to \infty \text{ as } x \to \infty \quad (12)$$

is sufficient. It should be observed that there are other possible sources of sequences of integers that satisfy (11). Specifically, it is known that the angles of the Gaussian integers $n = a^2 + b^2 = re^{i\theta}$, $\theta \in [0, \pi]$, are uniformly distributed for almost every integer $n$, see [ES99]. An infinite sequence of smooth or nearly powerful Gaussian integers with $\theta \approx \pi/2$ will be of interest here.

**Numerical Experiment I**

To test the density of the angles of the subsequence $q = p^{n+3} \leq x$ numerically, write $b_q = 2\sqrt{q} \cos\theta = q - \| q^{1/2} \|^2$, so that $\theta = \arccos\left( \frac{q - \| q^{1/2} \|^2}{2\sqrt{q}} \right) \in [0, \pi]$. Now following the idea of [LR70], the density of the proportion for which $a \leq \theta \leq b$ is given by

$$\lim_{x \to \infty} \frac{1}{x} \#\{ q \leq x : a \leq \theta \leq b \} = \int_a^b d\mu, \quad (13)$$

assuming that the probability measure μ exists. More advanced and current details on this topic are given in [IK04, p.487] and [MO06, p. 206]. The numerical experiment was geared to count the frequency of the angles in the interval [0, π]. The preliminary numerical data I have collected so far suggests that the angles are well distributed in [0, π] as expected and seem to be independent of the prime $p$, see (11).

**Numerical Experiment II**

The sequences of numbers $\{ N_n = q^n + 1 - a_n : n \in \mathbb{N} \}$ have rich mostly unexplored properties. The best known property is its divisibility: $N_d$ divides $N_n$ if $d$ divides $n$. The multiplicative structure is parametized by several parameters, for example, $a_n$, $q$, $n$. Thus, there are lots of possibilities of obtaining triples $a$, $b$, $c$ with large indices. And it is reasonable to expect one or more triples of large index $\gamma(a, b) > 1$ per prime powers.

The recursive formula $N_n = \#E(\mathbf{F}_{q^n}) = q^n + 1 - a_n$, where $a_0 = 2$, $a_1 = a_1$, and $a_n = a_1 a_{n-1} - q a_{n-2}$, $n \geq 2$, was used to conduct a light (not exhaustive) numerical experiment. It quickly produces many triples $a$, $b$, $c$ of indices $\gamma(a, b) > 1$, see the table below.
The indices of supersingular elliptic curves varies about $\gamma(a, b) = 1$, so these elliptic curves and nearly supersingular elliptic curves are excellent sources of smooth triples $a$, $b$, $c$.

| $a + b = c$ | Index $\gamma(a, b)$ | Genus $g$ |
|---|---|---|
| $2^8 + 1 - 4 = 5^3 \cdot$ | 1.426 | 1 |
| $2^{10} + 1 + 4 = 3 \cdot 7^3$ | 1.296 | 1 |
| $2^{11} + 1 + 138 = 3^7$ | 1.33 | 2 |
| $2^{12} + 1 + 102 = 3 \cdot 11^3$ | 0.94 | 2 |
| $2^{14} + 1 + 422 = 7^5$ | 1.117 | 2 |
| $2^{15} + 1 + 36 = 3^8 \, 5$ | 1.482 | 1 |
| $2^{16} + 1 + 88 = 3 \cdot 5^5 \cdot 7$ | 1.27 | 1 |
| $8^6 + 5^2 \cdot 41 = 3^6 \cdot 19^2$ | 1.240709 | 1 |
| $3^{14} + 5^4 \cdot 7 = 2^4 \cdot 547^2$ | 1.320123 | 1 |
| $3^{20} - 7 \cdot 16871 = 2^6 \cdot 11^4 \cdot 61^2$ | 1.099722 | 1 |
| $5^4 - 7^2 = 2^6 \cdot 3^2$ | 1.1887 | 1 |
| $7^8 - 4801 = 2^{10} \cdot 3^2 \, 5^4$ | 1.126071 | 1 |
| $11^6 + 2663 = 2^4 \cdot 3^4 \cdot 37^2$ | 0.917203 | 1 |
| $13^8 - 239^2 = 2^8 \cdot 3^2 \cdot 5^2 \cdot 7^2 \cdot 17^2$ | 1.26494 | 1 |
| $17^3 + 2^7 = 71^2$ | 1.094533 | 1 |
| $17^4 + 577 = 2^{10} \cdot 3^4$ | 1.03124 | 1 |
| $17^8 - 7^3 \cdot 487 = 2^{12} \cdot 3^4 \cdot 5^2 \cdot 29^2$ | 1.277956 | 1 |
| $19^6 + 3 \cdot 17 \cdot 269 = 2^4 \cdot 5^2 \cdot 7^6$ | 1.056668 | 1 |
| $23^6 + 5 \cdot 31 \cdot 157 = 2^6 \cdot 3^4 \cdot 13^4$ | 1.069422 | 1 |

## Appendix

This section reviews some results from the theory of integers in short intervals, and the related concept of powerful and smooth integers for the reader convenience, but no attempt was made to have a complete picture of the current state of knowledge. The results are stated verbatim (as in the original) for the reader convenience.

**Smooth Integers in Short Intervals**

For $x, y \in \mathbb{R}$, let $S(x, y) = \{ n \in \mathbb{N} : n \leq x \text{ and } p \mid n \Rightarrow p \leq y \}$ be the set of $y$-smooth integers $n \leq x$, and let $\psi(x, y) = \#S(x, y)$ be its cardinality.

***Theorem 5.*** ([BA01]) For all $x > x_0$, the interval $[\, x - x^{.525}, x \,]$ contains prime numbers.

***Theorem 6.*** ([HN91]) For any $\varepsilon > 0$ and $X > X_1(\varepsilon)$ the interval $[\, X, X + X^{1/2+\varepsilon} \,]$ contains an integer having no prime factor exceeding $e^{(\log X)^{2/3}+\varepsilon}$.

This improves a result of Balog which calls for prime factors $p \leq X^\varepsilon$. The improvement reduces the prime factors to subexponential size $p \leq e^{(\log X)^{2/3}+\varepsilon}$.

***Theorem 7.*** ([HN99]) Let $\varepsilon > 0$. Then for all large $x$, $\psi(x+x^{1/2},u) - \psi(x,u) \gg x^{1/2}$ provided that $u \geq x^{a/4+\varepsilon}$, $a = e^{-1/2} \approx .606$.

An integer $N$ is say to be $k$-powerful if $p^k$ divides $N$ for every prime factor $p$ of $N$.

***Theorem 8.*** ([HB87]) Every integer is a sum of at most three $k$-full integers.

***Theorem 9.*** ([DK05]) For any integer $k \geq 2$, there are infinitely many $N$, such that the open interval $(\, N^k, (N+1)^k \,)$ contains at least $M \geq [(3/8 + o(1))\log N / \log\log N]^{1/3}$ $k$-full integers.

***Theorem 10.*** ([LN02]) Suppose that $k$ is a finite field, and suppose that the cardinality of $k$ is an odd prime power $q$ and at least 14400. Then for all but at most $28q^{1/2}$ integers $z$ in the interval $[\, q^2 - 9^{-1}q^{3/2}, q^2 + 9^{-1}q^{3/2} \,]$ there are at least $q^{11/2}(48000\log(q)^2(\log\log(q)^2)$ sextic polynomials $f(x) \in k[x]$ with nonvanishing discriminants such that $\#J(k) = z$.